\documentstyle[leqno]{book}
\newcounter{supersection}[section]
\newtheorem{th}[supersection]{Theorem}
\newtheorem{lm}[supersection]{Lemma}
\newtheorem{re}[supersection]{Remark}
\newtheorem{co}[supersection]{Corollary}

\def\bibname{\textbf{REFERENCES}}
\def\thebibliography#1{\paragraph*{\uppercase{\bibname}}\list
{[\arabic{enumi}]}{\settowidth\labelwidth{[#1]}\leftmargin\labelwidth
\advance\leftmargin\labelsep\usecounter{enumi}}
\def\newblock{\hskip .11em plus .33em minus .07em}
\sloppy\clubpenalty4000\widowpenalty4000
\sfcode`\.=1000\relax}

\def\Dj{D{\hspace{-.75em}\raisebox{.3ex}{-}\hspace{.4em}}}

\begin{document}

\thispagestyle{plain}

\noindent {\small \sc $\;$}

\noindent {\small \sc $\;$}

\vspace*{10.00 mm}

\centerline{\Large \bf ONE METHOD FOR}

\medskip

\centerline{\Large \bf PROVING INEQUALITIES BY COMPUTER}
\footnotetext{Research partially supported by the MNTRS, Serbia, Grant No. 144020.}

\vspace*{6.00 mm}

\centerline{\large \it Branko J. Male\v sevi\' c}

\vspace*{4.00 mm}

\begin{center}
\parbox{25.0cc}{\scriptsize \bf
In this paper we consider a numerical method for proving a class of analytical inequalities via minimax
rational approximations. All numerical calculations in this paper are given by Maple computer program.}
\end{center}

\bigskip
\noindent
\section{\Large \bf \boldmath \hspace*{-7.0 mm}
1. Some particular inequalities} 

In this section we prove two new inequalities given in Theorem \ref{Th_first} and Theorem \ref{Th_second}.
While proving these theorems we use a method for inequalities of the following form $f(x) \geq 0$, for the
continues function \mbox{$f : [a,b] \longrightarrow R$}.

\bigskip
\noindent
\mbox{\large \textbf{1.1.}} 
Let us consider some inequalities for the gamma function which is defined by the integral:
\begin{equation}
\Gamma(z) = \displaystyle\int\limits_{0}^{\infty}{e^{-t} t^{z-1} \, dt}
\end{equation}
which converges for $Re(z) > 0$. In the paper \cite{Malesevic04} the following statement is proved.

\begin{lm}
For $x \in [0,1]$ the following inequalities are true$:$
\begin{equation}
\Gamma{\big (}x + 1{\big )}
<
x^2 - \displaystyle\frac{7}{4}x + \displaystyle\frac{9}{5}
\end{equation}
and

\medskip
\noindent
\begin{equation}
(x+2) \Gamma(x+1) > \displaystyle\frac{9}{5}.
\end{equation}
\end{lm}

\noindent
The previous statement (Lemma 4.1. of the paper \cite{Malesevic04}) is proved by the approximative
formula for the gamma function \mbox{$\Gamma(x+1)$} by the polynomial of the fifth order:
\begin{equation}
\quad
P_{5}(x)
=
-\mbox{\small $0.1010678$}x^5
+\mbox{\small $0.4245549$}x^4
-\mbox{\small $0.6998588$}x^3
+\mbox{\small $0.9512363$}x^2
-\mbox{\small $0.5748646$}x
+\mbox{\small $1$}
\end{equation}
which has the numerical bound of the absolute error $\varepsilon = 5 \cdot 10^{-5}$ for values of argument
$x \!\in\! [0,1]$ \cite{AbramowitzStegun72} (formula 6.1.35., page 257.).

\break

\medskip
\noindent
In the Maple computer program we use {\sf numapprox} package \cite{Geddes93} for obtaining
the minimax rational approximation $R(x) = P_{m}(x)/Q_{n}(x)$ of the continuous function $f(x)$
over segment $[a,b]$ {\big (}$m$ is the degree of the polynomial $P_{m}(x)$ and $n$ is the degree of
the polynomial $Q_{n}(x)${\big )}. Let $\varepsilon(x) = f(x) - R(x)$ be the error function of
an approximation over segment $[a,b]$. Numerical computation of $R(x)$ given by Maple command:

\smallskip

\noindent
\begin{equation}
R := \mbox{ \rm minimax $\!${\big (}$\!$ $f(x)$, $x=a$..$b$, $[m,n]$, '$err$' $\!${\big )} };
\end{equation}
The result of the previous command is the minimax rational approximation $R(x)$ and an estimate for
the value of the minimax norm of $\varepsilon(x)$ as the number $\mbox{\rm $err$}$ (computation is
realized without the weight function). With the Maple minimax command a realization of the Remez
algorithm is given \cite{W_Fraser_J_F_Hart}, \cite{C_Fike}. If it is not possible to determine
minimax approximation in Maple program there appears a message that it is necessary to increase
decimal degrees.

\medskip
\noindent
Let us assume that for the function $f(x)$ the minimax rational approximation $R(x)$ is determined.
Then, in the Maple the same estimate for the minimax norm of the error function
$\varepsilon(x)$ is given by command:
\begin{equation}
err:=\mbox{ \rm infnorm $\!{\big (}\!$ $\varepsilon(x)$, $x=a$..$b$ $\!${\big )} };
\end{equation}
The result of the previous command is number $err : = \mathop{\max}_{x \in [a,b]}{|f(x) - R(x)|}$.
Practically~for the bound of the absolute error function $|\varepsilon(x)|$ we use \mbox{$\varepsilon = err$}.
Let~us~remark that the bound of the absolute error $\varepsilon$ is a numerical bound in the sense
\cite{F_de_Dinechin_C_Lauter_G_Melquiond_05} (approximation errors on page 4.),
see also \cite{N_Brisebarre_JM_Muller_A_Tisserand_06}.

\medskip
\noindent
Let us notice, as it is emphasized by the Remark 4.2 of the paper \cite{Malesevic04}, that for the proof of Lemma 1.1.
it is possible to use other polynomial approximations (of~lower~degree) of the functions $\Gamma(x+1/2)$ and $\Gamma(x+1)$
for values $x \in [0,1]$. That idea is implemented in the next statement for the Kurepa's function which is defined
by the integral:
\begin{equation}
K(z)
=
\displaystyle\int\limits_{0}^{\infty}{
e^{-t} \displaystyle\frac{t^{z}-1}{t-1} \: dt},
\end{equation}
which converges for $Re(z) > 0$~\cite{Kurepa73}.
It is possible to make an analytical continuation of the Kurepa's function $K(z)$ to the meromorphic function
with simple poles at $z = -1$ and $z = -n$ $(n \geq 3)$. Practically for computation values of
the Kurepa's function we use the following formula:
\begin{equation}
K(z)
=
\displaystyle\frac{\mbox{\rm Ei}(1)+i\pi}{e}
+
\displaystyle\frac{(-1)^{z} \Gamma(1+z) \Gamma(-z,-1)}{e}
\end{equation}
which is cited in \cite{Malesevic03}. In the previous formula $\mbox{\rm Ei}(z)$ and $\Gamma(z,a)$ are
the exponential integral and the incomplete gamma function, respectively. Let us numerical prove
the following statement:

\begin{th}
\label{Th_first}
For $x \in [0,1]$ the following inequality is true$:$
\begin{equation}
K(x) \leq \displaystyle K^{'}\!(0) \, x ,
\end{equation}
where $K^{'}\!(0) = 1.432 \, 205 \, 735 \, \ldots$ is the best possible constant.
\end{th}

\noindent
{\bf Proof.}
Let us define the function $f(x) = K^{'}\!(0) \, x - K(x)$ for $x \in [0,1]$. Let us prove
$f(x) \geq 0$ for $x \in [0,1]$. Let us consider the continuous function:
\begin{equation}
g(x)
=
\left\{
\begin{array}{ccc}
\alpha                                                                    & : & x = 0,      \\[2.0 ex]
\displaystyle\frac{f(x)}{x^2}                                             & : & x \in (0,1];
\end{array}
\right.
\end{equation}
for constant $\alpha = -\displaystyle\frac{K^{''}\!(0)}{2}$. Let us notice that the constant:
\begin{equation}
\quad
\alpha
=\!
-\displaystyle\frac{K^{''}\!(0)}{2}
\,=\!
\lim\limits_{x \rightarrow 0+}{
\!\!\displaystyle\frac{K^{'}\!(0) - K^{'}\!(x)}{2x}}
\,=\!
\lim\limits_{x \rightarrow 0+}{
\!\!\displaystyle\frac{K^{'}\!(0)  \,x - K(x)}{x^2}}
\,=\!
\lim\limits_{x \rightarrow 0+}{
\!\!\displaystyle\frac{f(x)}{x^2}.}
\end{equation}
is determined in sense that $g(x)$ is a continuous function over segment $[0,1]$.
The numerical value of that constant is:
\begin{equation}
\alpha
=
-\displaystyle\frac{1}{2}
\displaystyle\int\limits_{0}^{\infty}{e^{-t} \frac{\log^{2}t}{t-1} \, dt}
=
0.963 \, 321 \, 189 \ldots  \, (>0) \; .
\end{equation}
Using Maple we determine the minimax rational approximation for the function $g(x)$ by the polynomial of the first order:
\begin{equation}
\quad
P_{1}(x) = -\mbox{\small $0.531 \, 115 \, 454$} \, x + \mbox{\small $0.921 \, 004 \, 887$}
\end{equation}
which has the bound of the absolute error $\varepsilon_{1} = 0.04232$ for values $x \!\in\! [0,1]$. The following
is true:
\begin{equation}
g(x) - {\big (} P_1(x) - \varepsilon_{1} {\big )} \geq 0
\quad \mbox{and} \quad
P_1(x) - \varepsilon_{1} > 0,
\end{equation}
for values $x \in [0,1]$. Hence, for $x \in [0,1]$ it is true that $g(x) > 0$, and $f(x) \geq 0$ as well.
\mbox{Q.E.D.}

\begin{re}
Numerical values of constants $K^{'}\!(0)$ and $K^{''}\!(0)$ are determined by Maple program.
The numerical value of $K^{'}\!(0)$ was first determined by D. Slavi\' c in~{\rm \cite{Slavic_73}}.
\end{re}

\begin{co}
A.$\,$Petojevi\' c in {\rm \cite{Petojevic05}} used an auxiliary result \mbox{$K(x) \!\leq\! 9/5 \, x$},
for values \mbox{$x \!\in\! [0,1]$}, from {\rm \cite{Malesevic04}} $($Lemma {\rm 4.3.}$)$, for proving
new inequalities for the Kurepa's function. Based on the previous theorem, all appropriate inequalities
from {\rm \cite{Petojevic05}} can be improved with a simple change of fraction $9/5$ with constant $K^{'}\!(0)$.
\end{co}

\bigskip
\noindent
\mbox{\large \textbf{1.2.}} 
D.$\,$S. Mitrinovi\' c considered in \cite{MitrinovicVasic70} the lower bound of the arc$\,$sin function,
which belongs to R.$\,$E. Shafer. Namely, the following statement is true.

\begin{th}
For $0 \leq x \leq 1$ the following inequalities are true$:$
\begin{equation}
\displaystyle\frac{3x}{2 + \sqrt{1-x^2}}
\leq
\displaystyle\frac{6(\sqrt{1+x} - \sqrt{1-x})}{4 + \sqrt{1+x} + \sqrt{1-x}}
\leq
\mbox{\rm arc$\,$sin} \, x \, .
\end{equation}
\end{th}

\medskip
\noindent
A.$\,$M. Fink proved the following statement in paper \cite{Fink95}.

\begin{th}
For $0 \leq x \leq 1$ the following inequalities are true$:$
\begin{equation}
\label{Ineq_Fink95}
\displaystyle\frac{3x}{2 + \sqrt{1-x^2}}
\leq
\mbox{\rm arc$\,$sin} \, x
\leq
\displaystyle\frac{\pi x}{2 + \sqrt{1-x^2}}.
\end{equation}
\end{th}

\medskip
\noindent
B.$\,$J. Male\v sevi\' c proved the following statement in \cite{Malesevic97}.

\begin{th}
For $0 \leq x \leq 1$  the following inequalities are true$:$
\begin{equation}
\label{Ineq_Malesevic97}
\displaystyle\frac{3x}{2 + \sqrt{1-x^2}}
\leq
\mbox{\rm arc$\,$sin} \, x
\leq
\displaystyle\frac{\mbox{\small $\displaystyle\frac{\pi}{\pi-2}$} x}{
\mbox{\small $\displaystyle\frac{2}{\pi-2}$} + \sqrt{1-x^2}}
\leq
\displaystyle\frac{\pi x}{2 + \sqrt{1-x^2}}.
\end{equation}
\end{th}

\begin{re}
The upper bound of the \mbox{arc$\,$sin} - function$:$
\begin{equation}
\phi(x)
=
\displaystyle\frac{\mbox{\small $\displaystyle\frac{\pi}{\pi-2}$} x}{
\mbox{\small $\displaystyle\frac{2}{\pi-2}$} + \sqrt{1-x^2}}
\end{equation}
is determined in paper {\rm \cite{Malesevic97}} by $\lambda$-method
Mitrinovi\' c-Vasi\' c {\rm \cite{MitrinovicVasic70}}.
\end{re}

\medskip
\noindent
L. Zhu proved the following statement in \cite{Zhu05}.

\begin{th}
For $x \in [0,1]$  the following inequalities are true$:$
\begin{equation}
\label{Ineq_Zhu05}
\begin{array}{rcl}
\displaystyle\frac{3x}{2 + \sqrt{1-x^2}}
&\!\!\!\leq\!\!\!&
\displaystyle\frac{6(\sqrt{1+x} - \sqrt{1-x})}{4 + \sqrt{1+x} + \sqrt{1-x}}
\,\leq\,
\mbox{\rm arc$\,$sin} \, x                                                       \\[3.0 ex]
&\!\!\!\leq\!\!\!&
\displaystyle\frac{\mbox{\small $\pi(\sqrt{2}+\displaystyle\frac{1}{2})$}
(\sqrt{1+x} - \sqrt{1-x})}{4 + \sqrt{1+x} + \sqrt{1-x}}
\,\leq\,
\displaystyle\frac{\pi x}{2 + \sqrt{1-x^2}}.
\end{array}
\end{equation}
\end{th}

\medskip
\noindent
In this paper we give an improved statement of L. Zhu.
Let us numerical prove the following statement:

\begin{th}
\label{Th_second}
For $x \in [0,1]$ the following inequalities are true$:$
\begin{equation}
\label{Ineq_Malesevic06}
\begin{array}{rcl}
\displaystyle\frac{3x}{2 + \sqrt{1-x^2}}
&\!\!\!\leq\!\!\!&
\displaystyle\frac{6(\sqrt{1+x} - \sqrt{1-x})}{4 + \sqrt{1+x} + \sqrt{1-x}}
\,\leq\,
\mbox{\rm arc$\,$sin} \, x                                                       \\[3.0 ex]
&\!\!\!\leq\!\!\!&
\displaystyle\frac{\mbox{\small $\displaystyle\frac{\pi(2-\sqrt{2})}{\pi-2\sqrt{2}}$}( \sqrt{1+x} - \sqrt{1-x} )}{
\mbox{\small $\displaystyle\frac{\sqrt{2}(4-\pi)}{\pi-2\sqrt{2}}$} + \sqrt{1+x} + \sqrt{1-x}}
                                                                                 \\[5.0 ex]
&\!\!\!\leq\!\!\!&
\displaystyle\frac{\mbox{\small $\pi(\sqrt{2}+\displaystyle\frac{1}{2})$}
(\sqrt{1+x} - \sqrt{1-x})}{4 + \sqrt{1+x} + \sqrt{1-x}}
\,\leq\,
\displaystyle\frac{\pi x}{2 + \sqrt{1-x^2}}.
\end{array}
\end{equation}
\end{th}

\noindent
{\bf Proof.} Inequality:
\begin{equation}
\displaystyle\frac{\mbox{\small $\pi(\sqrt{2}+\displaystyle\frac{1}{2})$}
(\sqrt{1+x} - \sqrt{1-x})}{4 + \sqrt{1+x} + \sqrt{1-x}}
\; \geq \;
\displaystyle\frac{\mbox{\small $\displaystyle\frac{\pi(2-\sqrt{2})}{\pi-2\sqrt{2}}$}( \sqrt{1+x} - \sqrt{1-x} )}{
\mbox{\small $\displaystyle\frac{\sqrt{2}(4-\pi)}{\pi-2\sqrt{2}}$} + \sqrt{1+x} + \sqrt{1-x}},
\end{equation}
for $x \in [0,1]$, is directly verifiable by algebraic manipulations. Let us define the following function:
\begin{equation}
f(x)
=
\displaystyle\frac{\mbox{\small $\displaystyle\frac{\pi(2-\sqrt{2})}{\pi-2\sqrt{2}}$}( \sqrt{1+x} - \sqrt{1-x} )}{
\mbox{\small $\displaystyle\frac{\sqrt{2}(4-\pi)}{\pi-2\sqrt{2}}$} + \sqrt{1+x} + \sqrt{1-x}}
-
\mbox{arc$\,$sin}\,x,
\end{equation}
for $x \in [0,1]$. Let us prove $f(x) \geq 0$ for $x \in [0,1]$, ie. $f(\sin t) \geq 0$ for
$t \in {\big [}0,\displaystyle\frac{\pi}{2}{\big ]}$. Let us define the function:
\begin{equation}
g(t) =
\left\{
\begin{array}{ccc}
\alpha        &\!\!:\!\!&                       t = 0,                                             \\[3.0 ex]
\displaystyle\frac{f(\sin t)}{
t^{3}{\Big (}\displaystyle \frac{\pi}{2}-t{\Big )}}                                                &\!\!:\!\!&
                                                t \in {\big (}0,\displaystyle\frac{\pi}{2}{\big )}, \\[3.5 ex]
\beta         &\!\!:\!\!&                       t = \displaystyle\frac{\pi}{2};
\end{array}
\right.
\end{equation}
where $\alpha$ and $\beta$ are constants determined with limits:
\begin{equation}
\alpha
\; =
\lim\limits_{t \rightarrow 0+}{\displaystyle\frac{f(\sin t)}{t^{3}{\Big (} \displaystyle\frac{\pi}{2}-t {\Big )}}}
\; = \;
\mbox{\small $\displaystyle\frac{(4 + \sqrt{2})\pi - 12\sqrt{2}}{(24-12\sqrt{2})\pi^2}$}
\; > \;
0
\end{equation}
and

\bigskip
\noindent
\begin{equation}
\beta
\; =
\lim\limits_{t \rightarrow \pi/2-}{\displaystyle\frac{f(\sin t)}{t^{3}{\Big (} \displaystyle\frac{\pi}{2}-t {\Big )}}}
\; = \;
\mbox{\small $\displaystyle\frac{(16\sqrt{2}-16) + (8-4\sqrt{2})\pi - \sqrt{2}\pi^2}{(2\sqrt{2}-2)\pi^3}$}
\; > \;
0.
\end{equation}
The previously determined function $g(t)$ is continuous over ${\big [}0,\displaystyle\frac{\pi}{2}{\big ]}$.
Using Maple we determine the minimax rational approximation for the function $g(t)$ by the polynomial of the first order:
\begin{equation}
P_{1}(t) = \mbox{\small $0.000 \, 410 \, 754$} \, t + \mbox{\small $0.000 \, 543 \, 606$}
\end{equation}
which has the bound of the absolute error $\varepsilon_{1} = 1.408 \cdot 10^{-5}$ for values
\mbox{$t \!\in\! {\big [}0,\displaystyle\frac{\pi}{2}{\big ]}$}. It is true:
\begin{equation}
g(t) - {\big (} P_{1}(t) - \varepsilon_{1} {\big )} \geq 0
\quad \mbox{and} \quad
P_{1}(t) - \varepsilon_{1} > 0,
\end{equation}
for values \mbox{$t \!\in\! {\big [}0,\displaystyle\frac{\pi}{2}{\big ]}$}. Hence, for
\mbox{$t \!\in\! {\big [}0,\displaystyle\frac{\pi}{2}{\big ]}$} is true that $g(t) > 0$
and therefore $f(\sin t) \geq 0$ for \mbox{$t \!\in\! {\big [}0,\displaystyle\frac{\pi}{2}{\big ]}$}.
Finally $f(x) \geq 0$ for $x \in [0,1]$.~\mbox{Q.E.D.}

\begin{re}
The paper {\rm \cite{Malesevic06b}} considers the upper bound of the \mbox{arc$\,$sin} - function$:$
\begin{equation}
\varphi(x)
=
\displaystyle\frac{\mbox{\small $\displaystyle\frac{\pi(2-\sqrt{2})}{\pi-2\sqrt{2}}$}( \sqrt{1+x} - \sqrt{1-x} )}{
\mbox{\small $\displaystyle\frac{\sqrt{2}(4-\pi)}{\pi-2\sqrt{2}}$} + \sqrt{1+x} + \sqrt{1-x}}
\end{equation}
via $\lambda$-method Mitrinovi\' c-Vasi\' c {\rm \cite{MitrinovicVasic70}}.
\end{re}

\section{\Large \bf \boldmath \hspace*{-7.0 mm}
2. A numerical method for proving inequalities} 

In this section we expose a numerical method for proving inequalities in following form:
\begin{equation}
\label{Ineq_Start}
f(x) \geq 0
\end{equation}
for the continuous function $f : [a,b] \longrightarrow R$. Let us assume that $x=a$ is the root of the order $n$
and $x=b$ is the root of the order $m$ of the function $f(x)$ (if $x=a$ is not the root then we determine that $n=0$,
ie. if $x=b$ is not the root then we determine that $m=0$). The method is based on the first assumption that
there exist finite and nonzero limits:
\begin{equation}
\label{Limits}
\alpha = \lim\limits_{x \rightarrow a+}{\displaystyle\frac{f(x)}{(x-a)^n(b-x)^m}}
\quad \mbox{and} \quad
\beta  = \lim\limits_{x \rightarrow b-}{\displaystyle\frac{f(x)}{(x-a)^n(b-x)^m}}.
\end{equation}
If for the function $f(x)$ (over extended domain of $[a, b]$) at the point $x = a$ there is an approximation
of the function by Taylor polynomial of $n$-th order and at point $x = b$ there is an approximation
of the function by Taylor polynomial of $m$-th order, then:
\begin{equation}
\alpha
=
\displaystyle\frac{f^{(n)}(a)}{n! \, (b-a)^{m}}
\qquad \mbox{and} \qquad
\beta
=
(-1)^m
\displaystyle\frac{f^{(m)}(b)}{m! \, (b-a)^{n}}.
\end{equation}
Let us define the function:
\begin{equation}
g(x)
 =
g^{f}_{a,b}{(x)}
=
\left\{
\begin{array}{ccc}
\alpha                                   &\!\!:\!\!& x = a                     \\[2.0 ex]
\displaystyle\frac{f(x)}{(x-a)^n(b-x)^m} &\!\!:\!\!& x \in (a,b)               \\[2.0 ex]
\beta                                    &\!\!:\!\!& x = b
\end{array}
\right.
\end{equation}
which is continuous over segment $[a,b]$. For proving inequality (\ref{Ineq_Start}) we use the equivalence:

\vspace*{1.5 mm}

\noindent
\begin{equation}
\label{Equivalence}
f(x) \geq 0 \Longleftrightarrow g(x) \geq 0,
\end{equation}
which is true for all values $x \in [a,b]$. Thus if $\alpha < 0$ or $\beta < 0$ the inequality (\ref{Ineq_Start})
is not true. Hence, we consider only the case $\alpha > 0$ and $\beta > 0$. Let us notice that if the function
$f(x)$ has only roots at some end-points of the segment $[a,b]$, then (\ref{Equivalence}) becomes $f(x) \geq 0$
iff $g(x) > 0$ for $x \in [a,b]$. The second assumption of the method is that there is the minimax (polynomial) rational
approximation $R(x) = P_{m}(x)/Q_{n}(x)$, of the function $g(x)$ over $[a,b]$, which has the bound of the absolute error
$\varepsilon > 0$ such that:

\vspace*{1.5 mm}

\noindent
\begin{equation}
\label{Ineq_Stop}
R(x) - \varepsilon > 0,
\end{equation}
for $x \in [a, b]$. Then $g(x) > 0$, for $x \in [a,b]$. Finally, on the basis (\ref{Equivalence}), we can conclude
that $f(x) \geq 0$, for $x \in [a,b]$.

\medskip
\noindent
Let us emphasize that the minimax (polynomial) rational approximation of the function $g(x)$ over $[a,b]$,
can be computed by Remez algorithm (via Maple minimax function \cite{Geddes93}). For applying Remez algorithm
to the function $g(x)$ it is sufficient that the function is continuous. If $g(x)$ is differentiable function
than the second Remez algorithm is applicable \cite{GolubSmith71}. According to the previous consideration,
the problem of proving inequality (\ref{Ineq_Start}), in some cases, becomes a problem of existence of the minimax
(polynomial) rational approximation $R(x)$ for $g = g^{f}_{a,b}(x)$ function with the bound of the absolute
error $\varepsilon > 0$ such that (\ref{Ineq_Stop}) is true. Let us notice that the problem of verification
of inequality (\ref{Ineq_Stop}) reduces to boolean combination of the polynomial inequalities.

\medskip
\noindent
Let us consider practical usages of the previously described method. For the function of one variable the previous
method can be applied to the inequality the infinite interval using the appropriate substitute variable, which
transforms inequality to the new one over the finite interval. Next, if some limits in (\ref{Limits}) are infinite,
then in some cases, the initial inequality can be transformed, by the means of appropriate substitute
variable, to the case when the both limits in (\ref{Limits}) are finite and nonzero.

\medskip
\noindent
The advantage of described method is that for the function $f(x)$ we don't have to use some regularities
concerning derivatives. Besides, the present method enables us to obtain computer-assisted proofs
of appropriate inequalities, which have been published in Journals which consider these topics.
With this method we have obtained numerical proofs of the appropriate inequalities
from the following articles
\cite{Elbert_Laforgia},
\cite{Cerone},
\cite{Guo_Qiao_Qi_Li}
\cite{Alzer},
\cite{Batir},
\cite{Chen_Qi},
\cite{Zhu},
\cite{Wu_Debentah},
\cite{SL_Qiu_MK_Vamanamurthy_M_Vuorinen},
\cite{SS_Dragomir_RP_Agarwal_NS_Barnett},
\cite{A_Laforgia_P_Natalini}.

\break

\medskip
\noindent
Finally, let us emphasize that the mentioned method can be extended and applied
to inequalities for multivariate functions by the means of
appropriate multivariate minimax rational approximations.

\bigskip
{\small
\noindent University of Belgrade,
          \hfill (Received   08/31/2006)                              \break
\noindent Faculty of Electrical Engineering,
          \hfill (Revised 10/30/2006)$\,$                             \break
\noindent P.O.Box 35-54, $11120$ Belgrade, Serbia               \hfill\break
\noindent {\footnotesize {\bf malesh@EUnet.yu}, {\bf malesevic@etf.bg.ac.yu}}
\hfill}

\end{document}